\documentclass[11pt,bezier]{article}
\usepackage{amsmath,amssymb,amsfonts,euscript}
\usepackage[usenames]{color}
\newcommand{\be}{\begin{equation}}
\newcommand{\ee}{\end{equation}}
\newtheorem{theorem}{Theorem}[section]

\newtheorem{corollary}{Corollary}[section]
\newtheorem{example}{Example}[section]
\newtheorem{definition}{Definition}[section]
\newtheorem{remark}{Remark}[section]

\renewcommand{\theequation}{\arabic{section}.\arabic{subsection}.\arabic{equation}}
\title{\bf\Large Characterization of Discrete Scale Invariant \\ Markov Sequences}
\vspace{-1cm}
\author
{ N. Modaresi\,\,\,\, and \,\,\,S. Rezakhah  \thanks{ Faculty of
Mathematics and Computer Science, Amirkabir University of
Technology, 424 Hafez Avenue, Tehran 15914, Iran. E-mail: $\;$
namomath@aut.ac.ir (N. Modaresi),$\;\;$ rezakhah@aut.ac.ir (S.
Rezakhah).} }
\date{}
\begin{document}
\maketitle

\begin{abstract}
By considering special sampling of discrete scale invariant (DSI) processes we provide a sequence which is in correspondence
to multi-dimensional self-similar process.
By imposing Markov property we show that the covariance functions of such discrete scale invariant Markov (DSIM) sequences
are characterized by variance, and covariance of adjacent samples in the first scale interval.
We also provide a theoretical method for estimating spectral density matrix of corresponding multi-dimensional self-similar Markov process.
Some examples such as simple Brownian motion with drift and scale invariant autoregressive model of order one are presented and these
properties are investigated. By simulating DSIM sequences we provide visualization of their behavior and investigate these results.
Finally we present a new method to estimate Hurst parameter of DSI processes and show that it has much better performance than maximum likelihood method for simulated data. 
\end{abstract}

\section{Introduction}
The notion of scale invariance or self-similarity is used as a fundamental property to handle and interpret many natural phenomena,
like textures in geophysics, turbulence of fluids, data of network traffic and image processing, etc \cite {b1}. The idea is that a
function is scale invariant if it is identical to any of its rescaled functions, up to some suitable renormalization of its
amplitude.\\
Discrete scale invariance (DSI) is a property which requires invariance by dilation for certain preferred scaling factors \cite {s1}.
It is known that DSI leads to log-periodic corrections to scaling. Log-periodic oscillations have been used to predict price trends,
turbulent time series, multi-fractal measures and crashes on financial markets \cite{z1}.
Burnecki et al. \cite{b4} studied $\alpha$-stable and self-similar processes and established the uniqueness of such process by
implementing Lamperti transform. They also present a natural construction of two distinct $\alpha$-stable Ornstein-uhlenbeck
processes via this transformation.
Borgnat et. al. \cite{b3} have studied the property of DSI and its relation to periodically correlated (PC) by means of the
Lamperti transformation too.
  Gray and Zhang proposed geometrically sampling of non-stationary self-similar processes, and implied
it to study Euler-Cauchy processes \cite{g3}. Vidacs and Virtamo considered geometric
sequence of sampling points in order to study the scaling behavior of the fractional Brownian motion traffic with fewer
samples \cite{v1}. They considered such sampling to obtain Maximum Likelihood estimator of the Hurst parameter of such models.
Current authors \cite{m2} considered some special geometric sampling to provide a correspondence multi-dimensional discrete
time self-similar process, and provided spectral representation and spectral density function of such
processes.\\
Continuous time DSI process in real life usually have some scale which is not necessarily integer,  so our special geometric sampling scheme is to obtain a sequence with DSI property.
Current authors \cite{m3} improved the efficiency of DSI sequence and considered some flexible sampling of a continuous time DSI process. They presented the spectral representation and  spectral density of such sampled process and implied this method to determine the structure of S$\&$P500 data for some special period.\\
We consider a DSI process with some scale $l>1$ and get our samples at points $\alpha^k$, where $k\in {\Bbb Z}$, $l=\alpha^T$ and
$T$ is the number of samples in each scale.
One of the advantages of such sampling is to provide a multi-dimensional self-similar process in correspondence to the geometrically sampled DSI sequence.
This provides an appropriate platform to study discrete scale
invariant processes. In this paper we study discrete scale invariant Markov (DSIM) sequences which have DSI property and are Markov in the wide sense.

We investigate  the properties of such DSIM sequences and show that its covariance function is characterized   by  $2T$ elements $\{R_{j}^H(1),R_{j}^H(0), j=0, 1,\ldots, T-1\}$, where $R_j^H(k)$ is the covariance function of $j$th and $(j+k)$th element of this sequence.

 We also present a theoretical estimation method for the spectral density matrix of corresponding multi-dimensional self similar Markov process and present a new method for estimation of Hurst parameter of DSI processes and implement the method for simulated data. \\
The paper is organized as follows. In section 2, we present definitions and some preliminary properties of discrete time DSI and self-similar processes, and also  Markov processes in
the wide sense.
In section 3, we present a characterization theorem for the covariance function of the DSIM sequences by applying the variance and covariance functions of the adjacent samples in the first scale interval.
By introducing discrete scale invariant autoregressive sequence of order $p$, DSIAR(p), with time varying coefficients and discrete
time simple Brownian motion we justify this result. Using quasi Lamperti transform, we show that DSIAR(1) is the counterpart of a
periodically correlated autoregressive, PCAR(1), process and obtain covariance function of it in this section.
In section 4, we present a characterization theorem  which can be used to verify whether a process has DSIM property. We
also introduce  a new method for the estimation of the covariance function which can be used to simply verify such property. This method enables us to estimate spectral density matrix of corresponding $T$-dimensional self-similar Markov process too.
In section 5, we present simulation of simple Brownian motion with different scale and Hurst parameters to visualize the behavior of
such DSIM process. We investigate the characterization of the covariance function of DSIM sequences for the simulated data to verify such DSIM property. Finally we present a new method in this section for estimating Hurst parameters for DSI process and shows its efficiency by comparing 
its performance  with the maximum likelihood method of corresponding multi-dimensional self-similar process for
the simulated data.
\vspace{-1mm}
\renewcommand{\theequation}{\arabic{section}.\arabic{equation}}
\section{Theoretical framework}
\setcounter{equation}{0}
\vspace{-3mm}
In this section, we review the concepts of self-similar, DSI and wide sense self-similar processes in discrete time.
We also present some characterizations of Markov processes in the wide sense.

\subsection{Discrete time self-similar processes}
A process $\{Y(t),t\in {\Bbb R}\}$ is said to be stationary, if for any $\tau\in {\Bbb R}$
\be\{ Y(t+\tau), t\in {\Bbb R}\}\stackrel{d}{=}\{Y(t), t\in {\Bbb R}\}\ee
where $\stackrel{d}{=}$ is the equality of all finite-dimensional distributions.
If $(2.1)$ holds for some $\tau\in {\Bbb R}$, the process is said to be periodically correlated (PC). The smallest of such $\tau$ is
called period of the process. Also a process $\{X(t),t\in {\Bbb R^+}\}$ is said to be self-similar of index $H$, if for any
${\lambda}>0$
\be \{X(t),t\in {\Bbb R^+}\}\stackrel{d}{=}\{\lambda^{-H}X(\lambda t),t\in {\Bbb R^+}\}.\ee
The process is said to be DSI of index $H$ and scaling factor ${\lambda}_0>0$ if $(2.2)$ holds
for $\lambda=\lambda_0$.

\begin{definition}
A process $\{X(k),k\in {\check{T}}\}$ is called discrete time self-similar process with parameter space $\check{T}$,
any subset of distinct points of real line, if for any $k_1, k_2 \in \check{T}$
\be \{X(k_2)\}\stackrel{d}{=}(\frac{k_2}{k_1})^H\{X(k_1)\}.\ee
The process $X(\cdot)$ is called DSI sequence with scale $l>0$ and parameter space $\check{T}$, if for any
$k_1, k_2=lk_1 \in \check{T}$, $(2.3)$ holds, see \em{\cite{m2}}.
\end{definition}

\begin{remark}
If the process $\{X(t),t\in {\Bbb R^+}\}$ is {\em DSI} with scale $l$, then for any fixed $s>0$, $X(\cdot)$ with parameter space
$\check{T}=\{l^{k}s; k\in {\Bbb Z}\}$ is a discrete time self-similar process.
If the process $\{X(t),t\in {\Bbb R^+}\}$ is {\em DSI} with scale $l=\alpha^T$, for some $T\in {\Bbb N}$ and $\alpha>1$, then for any
fixed $s>0$, $X(\cdot)$ with parameter space $\check{T}=\{\alpha^{k}s; k\in {\Bbb Z}\}$ is a {\em DSI} sequence, with the same scale.
\end{remark}

%\begin{remark}
%If the process $\{X(t),t\in {\Bbb R^+}\}$ is {\em DSI} with scale $l=\alpha^T$ for fixed $T\in {\Bbb N}$ and $\alpha>1$, then by
%sampling of the process at points $\alpha^{k}, k\in {\Bbb Z}$, we have $X(\cdot)$ as a {\em DSI} sequence with parameter space
%$\check{T}=\{\alpha^{k}; k\in {\Bbb Z}\}$ with scale $l$.
%If we consider sampling of $X(\cdot)$ at points $\alpha^{nT+k}, n\in {\Bbb Z} \mbox{for fixed}\ k=0, 1,\ldots, T-1$, then
%$X(\cdot)$ is a discrete time self-similar process with parameter space $\check{T}=\{\alpha^{nT+k}; n\in {\Bbb Z}\}$, see \cite{m2}.
%\end{remark}

\noindent Based on the definition of wide sense self-similar process presented in \cite{n1}, we present the following definition.
\begin{definition}
A random process $\{X(t),t\in {\Bbb R^+}\}$ is said to be discrete time self-similar in the wide sense with parameter
space $\check{T}=\{ \lambda^ns, n\in {\Bbb Z}\}$ for any $\lambda>1$, fixed $s>0$ and index $H>0$, if  the followings are satisfied
for every $n, n_1, n_2\in {\Bbb Z}$\\\vspace{-3mm}

$(i)\,\,\ E[X^2(\lambda^ns)]<\infty$,

$(ii)\,\,E[X(\lambda^{n+n_1}s)]=\lambda^{nH}E[X(\lambda^{n_1}s)]$,

$(iii)\,\, E[X(\lambda^{n+n_1}s)X(\lambda^{n+n_2}s)]=\lambda^{2nH}E[X(\lambda^{n_1}s)X(\lambda^{n_2}s)]$.\\
If the above conditions hold for fixed $n$, then the process is called {\em DSI} sequence in the wide sense with scale
$l=\lambda^{n}$. Another version of this definition is presented in \em{\cite{m2}}.
\end{definition}

Through this paper we are dealt with wide sense self-similar and wide sense scale invariant process, and for simplicity we omit
the term "in the wide sense" hereafter.

\subsection{Markov processes in the wide sense}
Let $\{X(n),n\in {\Bbb Z}\}$ be a second order process of centered random variables, $E[X(n)]=0$ and $E[|X(n)|^2]<\infty$.
Following Doob \cite {d1}, a real valued second order process $\{X(n),n\in {\Bbb Z}\}$ is Markov in the wide sense if, whenever
$t_1<\ldots<t_n$,\\\vspace{-3mm}
$$\hat{E}\big[X(t_n)|X(t_1),\ldots,X(t_{n-1})\big]=\hat{E}\big[X(t_n)|X(t_{n-1})\big]$$
is satisfied with probability one, where $\hat{E}$ stands for the linear projection (minimum variance estimator). If the process is
Gaussian, then $\hat{E}$ is a version of the conditional expectation.\\
The following fact on the covariance of Markov processes in the wide sense, are essentially due to Doob \cite {d1}. The stochastic process
$\{X(n),n\in {\Bbb Z}\}$ is Markov if and only $R(n_1,n)R(n,n_2)=R(n,n)R(n_1,n_2)$, for $n_1\leqslant n\leqslant n_2$ where
$R(n_1,n_2):=E[X(n_1)X(n_2)]$ is the covariance function of $X(\cdot)$.
It follows that $R(n_1,n_2)=G(n_1)K(n_2), n_1\leqslant n_2$
for some functions $G$ and $K$.
Borisov \cite {b2} completed the circle even for continuous time processes, namely, let $R(t_1,t_2)$ be some function defined on
$\mathcal{T}\times\mathcal{T}$ and suppose that $R(t_1,t_2)\neq 0$ everywhere on $\mathcal{T}\times\mathcal{T}$, where $\mathcal{T}$
is an interval. Then for $R(t_1,t_2)$ to be the covariance function of a Gaussian Markov process with time space $\mathcal{T}$ it
is necessary and sufficient that
\be R(t_1,t_2)=G\big(\min(t_1,t_2)\big)K\big(\max(t_1,t_2)\big)\ee
where $G$ and $K$ are defined uniquely up to a constant multiple and the ratio $G/K$ is a positive nondecreasing function on
$\mathcal{T}$.\\
It should be noted that the Borisov result on Gaussian Markov processes can be easily derived in the discrete case for second
order Markov processes in the wide sense, by using Theorem 8.1 of Doob \cite {d1}.

\renewcommand{\theequation}{\arabic{section}.\arabic{equation}}
\section{Covariance structure of the DSIM sequence}
\setcounter{equation}{0}
Let $\{X(t), t\in {\Bbb R^+}\}$ be a zero mean DSIM process with scale $l$. If $l<1$, we reduce the time scale,
so that $l$ be  greater than $1$ in the new scale. Now we consider to have $T$ samples in each scale. Our sampling scheme is
to get samples at points $1, \alpha, \alpha^2,\ldots, \alpha^{T-1},\alpha^T, \alpha^{T+1},\ldots$ where $\alpha$ is obtained by equality $l=\alpha^T$.
In this section we study the structure of the covariance function of DSIM sequence $\{X(\alpha^n), n\in{\Bbb Z}\}$.
We present a closed formula for the covariance function of DSIM sequence in Theorem 3.1. Discrete scale invariant autoregressive sequence of finite order $p$, DSIAR(p), is introduced. The covariance function of DSIM sequence under some conditions is characterized in Theorem 3.2. Finally we justify our result by presenting two examples of DSIM sequences as discrete time simple Brownian motion and DSIAR(1).

\begin{theorem}
Let $\{X(\alpha^n),n\in {\Bbb Z}\}$ be a {\em DSIM} sequence with scale $l=\alpha^T$, $\alpha>1$, $T\in {\Bbb N}$, then covariance
function
\be R_n^H(\tau):=E[X(\alpha^{n+\tau})X(\alpha^n)]\ee
where $\tau\in{\Bbb Z}$, $k\in {\Bbb W}=\{0, 1, \ldots \}$, $n=0, 1, \ldots, T-1$ and $R_n^H(\tau)\neq0$ is of the form
\be R_n^H(kT+v)=[\tilde{h}(\alpha^{T-1})]^k\tilde{h}(\alpha^{v+n-1})[\tilde{h}(\alpha^{n-1})]^{-1}R_n^H(0)\ee
where $v=0,1,\ldots,T-1$,
\be\tilde{h}(\alpha^r)=\prod_{j=0}^{r}R_j^H(1)/R_j^H(0),\hspace{7mm}r\in{\Bbb Z}\ee
and $\tilde{h}(\alpha^{-1})=1.$ Also
$R_n^H(-kT+v)=\alpha^{-2kTH}R_{n+v}^H((k-1)T+T-v).$
\end{theorem}
{\bf Proof:} From the Markov property (2.4), for $\alpha>1$, $R_n^H(\tau)$ satisfies\\\vspace{-3mm}
\be R_n^H(\tau)=G(\alpha^n)K(\alpha^{n+\tau}),\hspace{1cm}\tau=0, 1, \ldots\ee
By substituting $\tau=0$ in the above relation $G(\alpha^n)=\frac{R_{n}^H(0)}{K(\alpha^n)}$ and for $\tau\in{\Bbb W}$ we have
\be R_{n}^H(\tau)=\frac{K(\alpha^{n+\tau})}{K(\alpha^n)}R_{n}^H(0).\ee
So $K(\alpha^{n})=\frac{R_{n-1}^H(1)}{R_{n-1}^H(0)}K(\alpha^{n-1})$. Thus by iteration
\vspace{-3mm}
\be K(\alpha^n)=K(1)\prod_{j=0}^{n-1}h(\alpha^j)\ee
where $h(\alpha^j)=R_{j}^H(1)/R_{j}^H(0).$ As $X(\cdot)$ is a DSI sequence with parameter space $\{\alpha^k, k\in {\Bbb Z}\}$
and scale $\alpha^T$, by (3.1)
\be\frac{R_{T+i}^H(1)}{R_{T+i}^H(0)}=\frac{R_{i}^H(1)}{R_{i}^H(0)},\hspace{1cm}i= 0, 1, \ldots\ee
Thus for $n=0,1,\ldots,T-1,\hspace{2mm}k\in {\Bbb W}$, by (3.3) and the convention $\tilde{h}(\alpha^{-1})=1$ we have
\be K(\alpha^{kT+n})=K(1)\prod_{j=0}^{kT+n-1}h(\alpha^j)=K(1)[\tilde{h}(\alpha^{T-1})]^k\tilde{h}(\alpha^{n-1}).\ee
So by (3.5) for $\tau=kT+v$ and (3.8)
\be R_n^H(kT+v)=\frac{K(\alpha^{n+kT+v})}{K(\alpha^n)}R_{n}^H(0)
=[\tilde{h}(\alpha^{T-1})]^k\tilde{h}(\alpha^{v+n-1})[\tilde{h}(\alpha^{n-1})]^{-1}R_n^H(0)\ee
for $k\in\Bbb W$, $\alpha>1$ and $n,v=0, 1,\ldots, T-1$.\\ \\
Also using (3.1) for $\tau=-kT+v$
\be R_n^H(-kT+v)=E[X(\alpha^{-kT+n+v})X(\alpha^n)]=\alpha^{-2kTH}E[X(\alpha^{n+v})X(\alpha^{kT+n})]\ee
$$=\alpha^{-2kTH}R_{n+v}^H(kT+v)=\alpha^{-2kTH}R_{n+v}^H((k-1)T+T-v).\square$$

\begin{remark}
It follows from Theorem $3.1$ and relations $(3.3)$ and $(3.7)$ that $R_n^H(\tau)$, \\$n=0, 1,\ldots, T-1$ is fully specified
by the values of $\{R_{j}^H(1),R_{j}^H(0),j=0, 1, \ldots, T-1\}.$
\end{remark}

\subsection{Scale invariant autoregressive process}
First we present periodically correlated autoregressive (PCAR) process of finite order, see \cite{s0}, \cite{s01}. Then by introducing
discrete scale invariant autoregressive (DSIAR) sequence of order $p$, we obtain the covariance function of such process.\\\vspace{-3mm}
A PC Process $\{Y(n),n\in {\Bbb Z}\}$ is called causal autoregressive of order $p$, PCAR($p$), if\\
\be Y(n)=\phi_1(n)Y(n-1)+\phi_2(n)Y(n-2)+\ldots+\phi_p(n)Y(n-p)+Z(n)\ee
where $\{Z(n)\}$ is a periodic white noise with zero mean, variance $\sigma^2_n$ and period $T$ and the coefficients are periodic
with period $T$, that is
$\phi_i(n)=\phi_i(n+T)$ for $i= 1, \ldots, p$. By causality, we have that $Y(n)$ is uncorrelated with $\{Z(m), m>n\}$.
The PCAR(1) process is
\be Y(n)=\phi(n)Y(n-1)+Z(n)\ee
where in (3.11), $p=1$ and $\phi_1(n):=\phi(n)$. Then $\phi(n)=R_n(-1)/R_{n-1}(0)$, where $R_n(\tau)=E[Y(n+\tau)Y(n)]$.\\
We remind that a stochastic process $\{Y(n), n\in {\Bbb Z}\}$ is said to be $p$-ple Markov if for $n\geqslant p$
$$P\Big(Y(n)<x(n)|Y(v)=x(v), v\leqslant n-1\Big)$$
\vspace{-5mm}
\be =P\Big(Y(n)<x(n)|Y(n-1)=x(n-1), \ldots, Y(n-p)=x(n-p)\Big).\ee
The memory of this process extends $p$ units of time into the past. For the case of $p=1$ such a process is an ordinary Markov
process. Note that in the Gaussian case, the PC $p$-ple Markov is equivalent to the PCAR($p$).\\\

\noindent
{\bf Quasi Lamperti transform}\\
Here we present an extension of Lamperti transform which we call quasi Lamperti that enable us to introduce flexible sampling of a continuous DSI process and provide a DSI sequence by the followings.
\begin{definition}
The quasi Lamperti transform with positive index $H$ and $\alpha>1$ denoted by ${\cal L}_{H,\alpha}$ operates on a random process
$\{Y(t),t\in {\Bbb R}\}$ as
\be {\cal L}_{H,\alpha}Y(t)=t^HY(\log_\alpha t)\ee
and inverse quasi Lamperti transform ${\cal L}^{-1}_{H,\alpha}$ acts on process $\{X(t), t\in {\Bbb R^+}\}$ as
$${\cal L}^{-1}_{H,\alpha}X(t)={\alpha}^{-tH}X(\alpha^t).$$

\end{definition}

\begin{corollary}
If $\{X(t), t\in {\Bbb R^+}\}$ is {\em DSI} with scale ${\alpha}^T$ and $T\in {\Bbb N}$, then
$Y(t)={\cal L}^{-1}_{H,\alpha}X(t)$ is its {\em PC} counterpart  with period $T$. Conversely if $\{Y(t), t\in {\Bbb R}\}$
is {\em PC} with period $T$ then $X(t)={\cal L}_{H,\alpha}Y(t)$ is {\em DSI} with scale ${\alpha}^T$.
\end{corollary}
Let $X$ be a DSI process with scale $l>1$ and index $H$. We define discrete
scale invariant autoregressive sequence of order $p$, DSIAR(p), with parameter space $\check{T}=\{\alpha^k, k\in{\Bbb W}\}$
and scale $l=\alpha^T$ as
\be X(\alpha^n)=\theta_1(\alpha^{n-1})X(\alpha^{n-1})+ \ldots+ \theta_p(\alpha^{n-p})X(\alpha^{n-p})+\tilde{Z}(\alpha^n)\ee
where $\tilde{Z}(\cdot)$ is a white noise and $\theta(\cdot)$ has scale invariant property with scale $\alpha^T$, that is
$\theta_i(\alpha^{T+n})=\theta_i(\alpha^n)$ for $i=1, \ldots, p$.
Then we have $X(\alpha^{n+T})=\alpha^{TH}X(\alpha^n)$.

\begin{remark}
If $X$ is {\em DSI} with scale $l=\alpha$ and $\theta_i(\alpha^{n+1})=\theta_i(\alpha^n)$ for $i=1, \ldots, p$ then
$(3.15)$ is a self-similar autoregressive process of order $p$.
\end{remark}
Let $Y$ be the PC counterpart of $X$ that by Corollary 3.1 we have
$$X(\alpha^n)={\cal L}_{H,\alpha}Y(\alpha^n)=\alpha^{nH}Y(n).$$
In addition we assume that $\tilde{Z}(\alpha^n)={\cal L}_{H,\alpha}Z(\alpha^n)=\alpha^{nH}Z(n)$ and
$\theta_i(\alpha^{n-1})X(\alpha^{n-1})\\={\cal L}_{H,\alpha}W_i(\alpha^{n-1})$ for $i=1, \ldots, p$, where
$W_i(\alpha^{n-1})=\eta_i(\alpha^{n-1})Y(\alpha^{n-1})$ and $\eta_i(n-1)=\theta_i(\alpha^{n-1})$. So
$Y(n)=\alpha^{-H}\eta_1(n-1)Y(n-1)+ \ldots+ \alpha^{-pH}\eta_p(n-p)Y(n-p)+Z(n)$,
where $\eta_i(n)=\eta_i(n+T)$ is a PCAR(p) counterpart of $X$. If in (3.15) $p=1$, $X$ is a DSIAR(1)
process with parameter space $\check{T}=\{\alpha^k, k\in{\Bbb W}\}$ and $\theta_1(\alpha^n):=\theta(\alpha^n)$ as
\be X(\alpha^n)=\theta(\alpha^{n-1})X(\alpha^{n-1})+\tilde{Z}(\alpha^n).\ee
Then $\theta(\alpha^{n-1})=\frac{R_{n}^H(-1)}{R_{n-1}^H(0)}=\frac{R_{n-1}^H(1)}{R_{n-1}^H(0)}$, where $R_n^H(\cdot)$ is defined by
(3.1).
So its covariance function for $k\in\Bbb N$ and $\nu=0, 1, \ldots, T-1$ can be written as
$$R_n^H(kT+v)=E[X(\alpha^{kT+v+n})X(\alpha^n)]=\theta(\alpha^{kT+v+n-1})R_n^H(kT+v-1).$$
By a recursive method one can easily verify that
$R_n^H(kT+v)=\prod_{i=1}^{kT+v}\theta(\alpha^{kT+v+n-i})R_n^H(0)$.
Since $\theta(\cdot)$ has scale invariant property, then
\be R_n^H(kT+v)=\prod_{i=1}^{kT+v}\theta(\alpha^{v+n-i})R_n^H(0)
=\big[\prod_{j=1}^{T}\theta(\alpha^j)\big]^k\prod_{j=n}^{v+n-1}\theta(\alpha^j)R_n^H(0).\ee

\subsection{Characterization of the process}
By the following theorem, the necessary conditions are given to show that,
the covariance function  $R_n^H(\tau)$, given by (3.2) characterizes a DSIM sequence.

\begin{theorem}
The function $R_n^H(\tau)$ in relation $(3.2)$ characterize the covariance function of a {\em DSIM} sequence, if these functions
satisfy the following condition
\be \big(R_j^H(1)\big)^2\leqslant R_j^H(0)R_{j+1}^H(0)\ee
for $ j=0, 1,\ldots, T-1$ and $R_T^H(0)=\alpha^{2TH}R_0^H(0)$.
\end{theorem}
{\bf Proof:} It is enough to show that every covariance function of the form (3.2) is the covariance function of a DSIM sequence.
For this, we show in (i) that the covariance function $R_n^H(\tau)$ has scale invariance property
with scale $l=\alpha^T$.
We also prove in (ii) that $R_n^H(\tau)$ is the covariance function of a wide sense Markov process, as it satisfies in the
relation $(2.4)$.\\
{\bf (i)} According to $(3.2)$ we have
$$R_{n+T}^H(kT+v)=[\tilde{h}(\alpha^{T-1})]^k\tilde{h}(\alpha^{v+n+T-1})[\tilde{h}(\alpha^{n+T-1})]^{-1}R_{n+T}^H(0)$$
where by $(3.3)$ and $(3.7)$
$$\frac{\tilde{h}(\alpha^{v+n+T-1})}{\tilde{h}(\alpha^{n+T-1})}
=\frac{\prod_{j=T}^{v+n+T-1}\frac{R_j^H(1)}{R_j^H(0)}}
{\prod_{j=T}^{n+T-1}\frac{R_j^H(1)}{R_j^H(0)}}
=\frac{\prod_{j=0}^{v+n-1}\frac{R_{j+T}^H(1)}{R_{j+T}^H(0)}}{
\prod_{j=0}^{n-1}\frac{R_{j+T}^H(1)}{R_{j+T}^H(0)}}
=\frac{\tilde{h}(\alpha^{v+n-1})}{\tilde{h}(\alpha^{n-1})}$$
and $R_{n+T}^H(0)=l^{2H}R_n^H(0)$, therefore $R_{n+T}^H(kT+v)=l^{2H}R_n^H(kT+v)$.\\
{\bf (ii)} By (3.1), (3.2) and (3.8) we have that
$$R_n^H(kT+v)=K(1)\tilde{h}(\alpha^{kT+n+v-1})[K(1)\tilde{h}(\alpha^{n-1})]^{-1}R_n^H(0)$$
$$=K(\alpha^{kT+n+v})[K(\alpha^n)]^{-1}R_n^H(0)=K(\alpha^{kT+n+v})G(\alpha^n)$$
where $G(\alpha^n)=[K(\alpha^n)]^{-1}R_n^H(0)$. Thus for $\tau\in {\Bbb Z}$, $\alpha>1$ we have
$R_n^H(\tau)=K(\alpha^{n+\tau})G(\alpha^n)$,
that satisfies Borisov condition (2.4). So $R_n^H(\cdot)$ is the covariance function of a Markov process, provided $G/K$ is
positive and nondecreasing. Positivity of $G/K$ is straight as
\be\frac{G(\alpha^n)}{K(\alpha^{n})}=\frac{R_n^H(0)}{K^2(\alpha^{n})}>0.\ee
Now we prove that $F(\alpha^r)=G(\alpha^r)/K(\alpha^r)$ under condition (3.18) is nondecreasing function. So we show that
$F(\alpha^{r+1})/F(\alpha^r)\geqslant 1$. Let $r=kT+n\geqslant0$ and $0\leqslant n\leqslant T-1$, then we consider two cases,
$0\leqslant n\leqslant T-2$ and $n=T-1$.
For $0\leqslant n\leqslant T-2$, by (3.8) and (3.19) we have
$$\frac{F(\alpha^{r+1})}{F(\alpha^r)}=\frac{K^2(\alpha^{kT+n})}{K^2(\alpha^{kT+n+1})}\frac{R_{kT+n+1}^H(0)}{R_{kT+n}^H(0)}=
\Big[\frac{\tilde{h}(\alpha^{n-1})}{\tilde{h}(\alpha^{n})}\Big]^2
\frac{R_{n+1}^H(0)}{R_{n}^H(0)}.$$
As $X(\cdot)$ is DSI sequence with parameter space $\{\alpha^k, k\in {\Bbb Z}\}$ and scale $\alpha^T$, so by (3.1),
$R_{kT+n}^H(0)=\alpha^{2kTH}R_n^H(0)$. Therefore by (3.3)
$$\frac{F(\alpha^{r+1})}{F(\alpha^r)}=\Big[\frac{R_{n}^H(0)}{R_n^H(1)}\Big]^2\frac{R_{n+1}^H(0)}{R_{n}^H(0)}
=\frac{R_{n}^H(0)R_{n+1}^H(0)}{(R_{n}^H(1))^2}.$$
Under condition $(3.18)$ for $j=n$, $F(\alpha^{r+1})/F(\alpha^r)\geqslant 1$, so $F(\alpha^r)$ is nondecreasing.
For $n=T-1$
$$\frac{F(\alpha^{r+1})}{F(\alpha^r)}=\frac{K^2(\alpha^{kT+T-1})}{K^2(\alpha^{kT+T})}\frac{R_{kT+T}^H(0)}{R_{kT+T-1}^H(0)}
=[\frac{\tilde{h}(\alpha^{T-2})}{\tilde{h}(\alpha^{T-1})}]^2\frac{R_{T}^H(0)}{R_{T-1}^H(0)}
=\frac{R_{T-1}^H(0)R_{T}^H(0)}{(R_{T-1}^H(1))^2}.$$
By a similar method and under condition $(3.10)$ for $j=T-1$, $F(\alpha^r)$ is nondecreasing.$\square$

\subsection{Examples}
We present two examples of DSIM sequence as discrete time simple Brownian motion with drift and DSIAR(1), then justify Theorem 3.1.

\begin{example}
A process $X(t)$ is simple Brownian motion, with drift, and index $H>0$, $a>0$ and scale $\lambda>1$ as
\vspace{-3mm}
$$X(t)=\sum_{n=1}^{\infty}\lambda^{n(H-\frac{1}{2})}I_{[\lambda^{n-1}, \lambda^{n})}(t)\big[B(t)+\lambda^{n/2}g(\lambda^{-n}t)\big]$$
where $B(\cdot)$, $I(\cdot)$ are Brownian motion and indicator function respectively and $g(\cdot)$ is a real deterministic function.
\end{example}
Let $A_n=[\lambda^{n-1}, \lambda^{n})$, $n\in{\Bbb N}$ be disjoint sets.
The expectation value and covariance function of the process for $t\in A_n$, $s\in A_m$ and $s\leqslant t$ are
$E(X(t))=\lambda^{n(H-\frac{1}{2})}\big[E(B(t))+\lambda^{n/2}g(\lambda^{-n}t)\big]=\lambda^{nH}g(\lambda^{-n}t),$
\vspace{-3mm}
$$\mathrm{Cov}\big(X(t), X(s)\big)=\lambda^{(n+m)H'}\mathrm{Cov}\big(B(t)+\lambda^{n/2}a,
B(s)+\lambda^{m/2}a\big)=\lambda^{(n+m)H'}s$$
where $\mathrm{Cov}\big(B(t), B(s)\big)=\min\{t,s\}$ and $H'=H-\frac{1}{2}$. Therefore by the condition $(2.4)$, the above covariance
is the covariance function of a Markov process. If  $t\in A_{n+1}$ and $s\in A_{m+1}$ we have $E(X(\lambda t))=\lambda^{H}E(X(t))$
and
\vspace{-1mm}
$$\mathrm{Cov}\big(X(\lambda t), X(\lambda s)\big)=\lambda^{(n+m+2)H'}\lambda s=\lambda^{2H}\mathrm{Cov}\big(X(t), X(s)\big).$$
Then $X(t)$ is DSI with scale $\lambda$.
By sampling of the process $X(\cdot)$ at points $\alpha^n$, $n\in {\Bbb W}$, where $\lambda=\alpha^T$, $T\in {\Bbb N}$ and
$\lambda>1$,
we provide a DSIM sequence and investigate the conditions of Theorem $3.1$.\\
For $j=kT+i$ where $i=0, 1, \ldots, T-2$ and $k\in\Bbb W$ we have
$$h(\alpha^j)=\frac{R_j^H(1)}{R_j^H(0)}=\frac{\mathrm{Cov}\big(X(\alpha^{j+1}),
X(\alpha^j)\big)}{\mathrm{Cov}\big(X(\alpha^j), X(\alpha^j)\big)}=\frac{\alpha^{2(k+1)TH'+j}}{\alpha^{2(k+1)TH'+j}}=1$$
as $\alpha^j, \alpha^{j+1} \in A_{k+1}$. Also for $j=kT+T-1$ when $\alpha^j\in A_{k+1}$ and $\alpha^{j+1}\in A_{k+2}$,
we have $h(\alpha^j)=\alpha^{TH'}$. Thus for $j=kT+i$, $i=0, 1, \ldots, T-2$ and $k\in\Bbb W$
$$\tilde{h}(\alpha^{kT+i})=\prod_{r=0}^{kT+i}h(\alpha^r)=\prod_{r=0}^{k}\alpha^{TH'}=\alpha^{kTH'}.$$
For $j=kT+T-1$, $\tilde{h}(\alpha^{kT+T-1})=\alpha^{(k+1)TH'}$,
$$\tilde{h}(\alpha^{v+n-1})=\left\{\begin {array}{cc}
1\hspace{2cm}v+n-1\leqslant T-2\\
\alpha^{TH'}\hspace{15mm}v+n-1\geqslant T-1\\
\end {array}\right.$$\\
and $\tilde{h}(\alpha^{n-1})=1$, $R_n^H(0)=E[X^2(\alpha^n)]=\alpha^{2TH'+n}$, $n=0, 1, \ldots, T-1$. Thus by (3.2)
$$R_n^H(kT+v)=\left\{\begin {array}{cc}
\alpha^{(k+2)TH'+n}\hspace{1cm}v+n-1\leqslant T-2\\
\alpha^{(k+3)TH'+n}\hspace{1cm}v+n-1\geqslant T-1\\
\end {array}\right.$$\\
Also by straight calculation we have the same result.

\begin{example}
We show that the covariance function of \em{DSIAR(1)} defined by $(3.16)$
satisfies the relation $(3.2)$. By using $(3.3)$ and by the fact that $\theta(\alpha^{n-1})=R_n^H(-1)/R_{n-1}^H(0)$ we have
$$h(\alpha^{n-1})=\frac{R_{n-1}^H(1)}{R_{n-1}^H(0)}=\frac{E[X(\alpha^n)X(\alpha^{n-1})]}{E[X(\alpha^{n-1})X(\alpha^{n-1})]}
=\frac{R_{n}^H(-1)}{R_{n-1}^H(0)}=\theta(\alpha^{n-1})$$
then according to $(3.2)$,
$R_n^H(kT+v)=[\tilde{\theta}(\alpha^{T-1})]^k\tilde{\theta}(\alpha^{v+n-1})[\tilde{\theta}(\alpha^{n-1})]^{-1}R_n^H(0)$
where $\tilde{\theta}(\alpha^r)=\prod_{j=1}^{r}\theta(\alpha^j)$.Therefore\vspace{-3mm}
$$R_n^H(kT+v)=\big[\prod_{j=0}^{T-1}\theta(\alpha^j)\big]^k\prod_{j=n}^{n+v-1}\theta(\alpha^j)R_n^H(0)$$
which is the same as the straight computation of the covariance function in $(3.17)$.
\end{example}

\section{Spectral density estimation}
\setcounter{equation}{0}
In this section we assume that $\{X(\alpha^n), n\in {\Bbb Z}\}$ is a DSIM sequence with scale $l=\alpha^T$ and present the characterization Theorem 4.1 for the covariance function of the associated T-dimensional discrete time self-similar Markov process. Then we report the spectral density matrix from \cite{m2}, and present a dynamic method for estimation of the covariance function of DSIM process and spectral density matrix of its corresponding multi-dimensional self-similar process.
The following definition, remark and theorem are reported from \cite{m2}.

\begin{definition}
The process $U(t)=(U^0(t), U^1(t),\ldots,U^{q-1}(t))$ with parameter space $\check{T}=\{l^{n}, n\in {\Bbb Z}\}$, $l=\alpha^T$,
$\alpha>1$ and $T\in {\Bbb N}$ is a q-dimensional discrete time self-similar process in the wide sense, where

$\Bbb (a)$\hspace{3mm} $\{U^j(\cdot)\}$ for all $j=0, 1, \ldots, q-1$ is discrete time self-similar process with parameter space
$\check{T}^j=\{l^{n}, n\in {\Bbb Z}\}$.
%\vspace{-2mm}

$\Bbb (b)$\hspace{3mm} For every $n, \tau\in {\Bbb Z},\,\ j, k=0, 1, \ldots, q-1$
$$\mathrm{Cov}\big(U^j(l^{n+\tau}), U^k(l^{n})\big)=l^{2nH}\mathrm{Cov}\big(U^j(l^{\tau}),U^k(1)\big).$$
\end{definition}

\begin{remark}
Corresponding to the {\em DSIM} sequence $\{X(\alpha^n),n\in {\Bbb Z}\}$ with scale $l=\alpha^T$, $\alpha>1$ there exists a
$T$-dimensional discrete time self-similar Markov process $V(t)=\big(V^0(t), V^1(t), \ldots,\\V^{T-1}(t)\big)$ with parameter
space $\check{T}=\{l^{n}, n\in {\Bbb Z}\}$ and
\be V^k(l^n)=V^k(\alpha^{nT}):=X(\alpha^{nT+k}),\hspace{7mm}k=0,1,\ldots,T-1\ee
\end{remark}
This remark is valid by the fact that $\{V^k(l^n), n\in{\Bbb Z}\}$ for $k=0, \ldots, T-1$ are discrete time
self-similar process, and for $n, n_1, n_2\in {\Bbb Z}$
$$\mathrm{Cov}\big(V^j(l^{n+n_1}),V^i(l^{n+n_2})\big)=\mathrm{Cov}\big(X(\alpha^{(n+n_1)T+j}),X(\alpha^{(n+n_2)T+i})\big)$$
\vspace{-8mm}
$$=\alpha^{2nTH}\mathrm{Cov}\big(X(\alpha^{n_1T+j}),X(\alpha^{n_2T+i})\big)=l^{2nH}\mathrm{Cov}\big(V^j(l^{n_1}),
V^i(l^{n_2})\big).$$
So assertions (a) and (b) of Definition 4.1 are satisfied.\\

Let $Q^H(n,\tau)=[Q^H_{jk}(n,\tau)]_{j,k=0, 1,\ldots, T-1}$ be the covariance
matrix of $V(l^n)$, then
$$Q^H_{jk}(n,\tau)=E[V^j(l^{n+\tau})V^k(l^n)]=E[X(\alpha^{(n+\tau)T+j})X(\alpha^{nT+k})].$$

\begin{theorem}
Let $\{X(\alpha^n),n\in {\Bbb Z}\}$ be a {\em DSIM} sequence with the covariance function $R_n^H(\tau)$ and let
$\{V(l^n),n\in {\Bbb Z}\}$, defined in $(4.1)$, be its associated T-dimensional discrete time self-similar process with
covariance function $Q^H(n,\tau)$. Then
\be Q^H(n,\tau)=\alpha^{2nHT}Q^H(\tau)=\alpha^{2nHT}C^HR^H[\tilde{h}(\alpha^{T-1})]^{\tau},\hspace{7mm}\tau\in{\Bbb Z}\ee
where $\tilde{h}(\cdot)$ is defined by $(3.3)$ and the matrix $C^H$ is given by $C^H=[C^H_{jk}]_{j,k= 0, \ldots, T-1}$,
where $C^H_{jk}=\tilde{h}(\alpha^{j-1})[\tilde{h}(\alpha^{k-1})]^{-1}$, and the matrix $R^H$ is a diagonal matrix with
diagonal elements $R^H_{j}(0)$, $j=0, 1, \ldots, T-1$.
\end{theorem}
%\red{{\bf Proof:}
%As $X(\cdot)$ is DSI sequence with parameter space $\{\alpha^k, k\in {\Bbb Z}\}$ and scale $\alpha^T$, so
%\be Q^H_{jk}(n,\tau)=\alpha^{2nHT}E[X(\alpha^{\tau T+j})X(\alpha^{k})]=\alpha^{2nHT}R_k^H(\tau T+j-k).\ee
%Then $Q^H_{jk}(n,\tau)=\alpha^{2nHT}Q^H_{jk}(\tau)$ and $Q^H(\tau)=[Q^H_{jk}]_{j,k= 0, \ldots, T-1}$ where
%$Q^H_{jk}(\tau)=R_k^H(\tau T+j-k)$.
%By the Markov property of $X(\cdot)$ from Theorem $3.1$, for
%$0\leqslant j-k\leqslant T-1$
%$$R_k^H(\tau T+j-k)=[\tilde{h}(\alpha^{T-1})]^{\tau}\tilde{h}(\alpha^{j-1})[\tilde{h}(\alpha^{k-1})]^{-1}R_k^H(0).$$
%Let $C^H_{jk}=\tilde{h}(\alpha^{j-1})[\tilde{h}(\alpha^{k-1})]^{-1}$, so
%\be Q^H_{jk}(\tau)=[\tilde{h}(\alpha^{T-1})]^{\tau}C^H_{jk}R_k^H(0),\hspace{1cm}\tau\in {\Bbb Z}.\square\ee}

Current authors, showed that \cite {m2} if
$\{X(\alpha^n), n\in{\Bbb Z}\}$ is a discrete time self-similar process with scale $l=\alpha^T$, $T\in {\Bbb N}$ then
the spectral representation of the covariance function of the process is\vspace{-4mm}
\be R_n^H(\tau):=\mathrm{Cov}\big(X(\alpha^n),X(\alpha^{n+\tau})\big)=\alpha^{(2n+\tau)H}\sum_{k=0}^{T-1}B_k(\tau)e^{2k\pi in/T}\ee
where
\vspace{-5mm}
\be B_k(\tau)=\int_{0}^{2\pi}e^{i\tau\omega}f_k(\omega)d\omega\ee
and
\be f_{jk}(\omega)=\frac{1}{T}f_{k-j}\big((\omega-2\pi j)/T\big)\ee
for $j,k=0,1,\ldots,T-1$ and $0\leqslant\omega<2\pi$.\\
Using Definition 4.1  they proved that the spectral density matrix of such $T$-dimensional process is
$d^H(\omega)=[d_{jr}^H(\omega)]_{j,r=0, 1, \ldots, T-1}$, where
\be d^H_{jr}(\omega)=\frac{1}{2\pi}\left[\frac{\tilde{h}(\alpha^{j-1})R_r^H(0)}{\tilde{h}(\alpha^{r-1})
(1-e^{ -i\omega T}\alpha^{- H T}\tilde{h}({\alpha}^{T-1}))}-
\frac{\tilde{h}(\alpha^{r-1})R_j^H(0)}{\tilde{h}(\alpha^{j-1})
\big(1-e^{-i\omega T}\alpha^{HT}\tilde{h}^{-1}({\alpha}^{T-1})\big )}\right],\ee
and $\tilde{h}(\cdot)$ is defined by (3.3). Thus we have the following result.

\begin{remark}
As the spectral density matrix of the multi-dimensional self-similar process $V(l^n)$, defined by $(4.1)$ is characterized by
$\{R_{j}^H(1),R_{j}^H(0), j=0, 1, \ldots, T-1\}$, relations $(4.3)-(4.5)$ reveal that the spectral density
of the corresponding {\em DSIM} sequence $\{X(\alpha^k), \; k\in {\Bbb Z}\}$ is fully specified by
$\{R_{j}^H(1),R_{j}^H(0), j=0, 1, \ldots, T-1\}$.
\end{remark}

\begin{example}
We present the T-dimensional {\em DSIM} sequence corresponding to the \\{\em DSIAR(1)}, defined by $(3.16)$ as
$V(l^n)=(V^0(l^n), V^1(l^n), \ldots, V^{T-1}(l^n))$, where $V^k(l^n)=X(\alpha^{nT+k})$. Also we have
$$h(\alpha^j)=\frac{R_j^H(1)}{R_j^H(0)}=\frac{E[X(\alpha^{j+1})X(\alpha^j)]}{E[X(\alpha^j)X(\alpha^j)]}
=\frac{R_{j+1}^H(-1)}{R_j^H(0)}=\theta(\alpha^j).$$
Thus the spectral density matrix of $V(l^n)$ is obtained by substituting of $\tilde{h}(\alpha^j)=\tilde{\theta}(\alpha^j)$ in
$(4.6)$.
\end{example}

\noindent
{\bf Estimation of covariance functions and spectral density matrix}\\
Here we explain our method for estimation of the spectral density matrix of multi-dimensional self-similar Markov processes presented by (4.6).
For this, we need to estimate $R_j^H(0)$ and
$R_j^H(1)$ for $j=0, \ldots, T-1$ and $R_n^H(0)$.
By the DSI property of the process, $X(\alpha^j)$ and $\lambda^{-kH}X(\alpha^{kT+j})$ have the same distribution,
so we evaluate these estimations by the followings. Let
$$m_j=\sum_{k=0}^{M-1}\frac{X(\alpha^{kT+j})}{\lambda^{k}M}.$$
Then we use following estimations to check the relation (3.2)
$$\hat{R}_j^H(0)=\widehat{\mbox{Var}}\big(X(\alpha^j)\big)=\sum_{k=0}^{M-1}\big(\lambda^{-kH} X(\alpha^{kT+j})- m_j\big)^2/(M-1)$$
for $j=0, \ldots, T-1$.
$$\hspace{-1.3in}\hat{R}_j^H(1)=\widehat{\mbox{Cov}}\big(X(\alpha^j), X(\alpha^{j+1}\big)=\hspace{1.3in}$$
$$\sum_{k=0}^{M-1}\big(\lambda^{-kH} X(\alpha^{kT+j})- m_j\big)\big(\lambda^{-kH} X(\alpha^{kT+j+1})- m_{j+1}\big)/(M-1)$$
for $j=0, \ldots, T-2$, and $\hat{R}_{T-1}^H(1)=\widehat{\mbox{Cov}}\big(X(\alpha^{T-1}), X(\alpha^{T})\big)$ where
$$\hat{R}_{T-1}^H(1)=\frac{1}{M-1}\sum_{k=0}^{M-1}
\big(\lambda^{-kH} X(\alpha^{kT+T-1})- m_{T-1}\big)\big(\lambda^{-kH} X(\alpha^{kT+T})-\lambda^Hm_{0}\big).$$
Also for $n=rT+i$ , $0\leq i \leq T-1$ we have that
$$\hat{R}_n^H(0)=\widehat{\mbox{Var}}\big(X(\alpha^n)\big)=\frac{1}{M-r-1}\sum_{k=0}^{M-r-1}
\big(\lambda^{-kH} X(\alpha^{kT+n})-\lambda^{rH}m_i\big)^2$$
and
$$\hspace{-1in} \hat{R}_n^H(\tau)=\widehat{\mbox{Cov}}\big(X(\alpha^n), X(\alpha^{n+\tau})\big)=\frac{1}{M-s-1}\hspace{1in}$$
$$\sum_{k=0}^{M-s-1}\big(\lambda^{-kH} X(\alpha^{kT+n})- \lambda^{rH}m_i\big)\big(\lambda^{-kH} X(\alpha^{kT+n+\tau})-
\lambda^{sH}m_j\big)$$
where $n+\tau=sT+j$ and  $0\leq j\leq T-1$. By applying the above estimators in (4.6), the estimation of spectral density matrix $d^H(\omega)$ is evaluated.

\section{Simulation and Estimation}
In this section first we present simulation of simple Brownian motion with different scale and Hurst parameters to
visualize the behavior of such DSIM process in subsection 5.1. We verify the main results of the paper as
Theorems 3.1 and 3.2, by simulating such process and estimating covariance and variance functions involve in relation (3.2).
We use the scale invariant property of the process, with known scale parameter in subsection 5.2. So this study would be
verification of Markov property of scale invariant processes. Finally in subsection 5.3, we present a new method for
estimating Hurst parameters of DSI and self-similar processes.

\input{epsf}
\epsfxsize=3in \epsfysize=1.5in
\begin{figure}
\centerline{ \epsffile{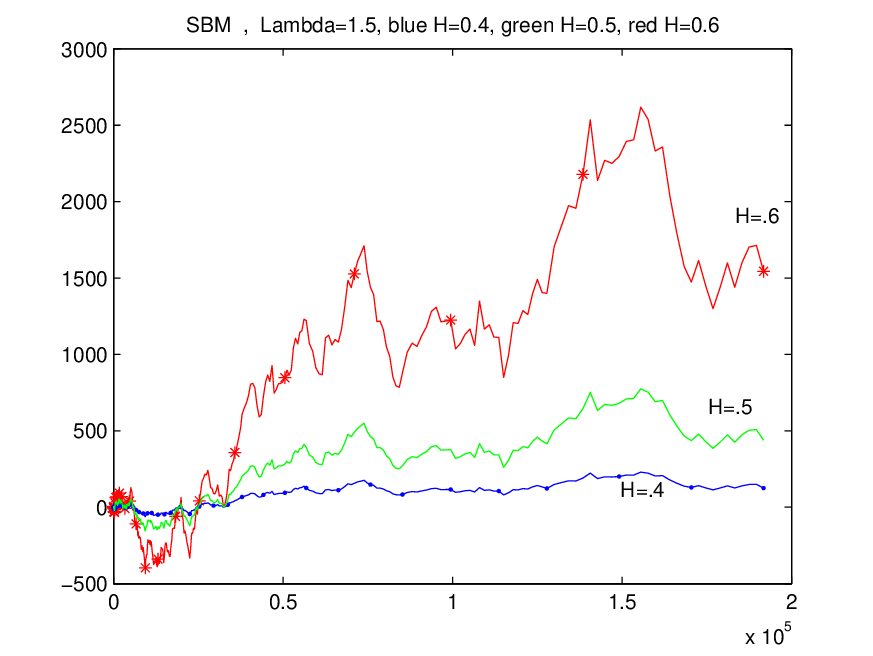}\hspace{.1in}\epsfxsize=3in \epsfysize=1.5in \epsffile{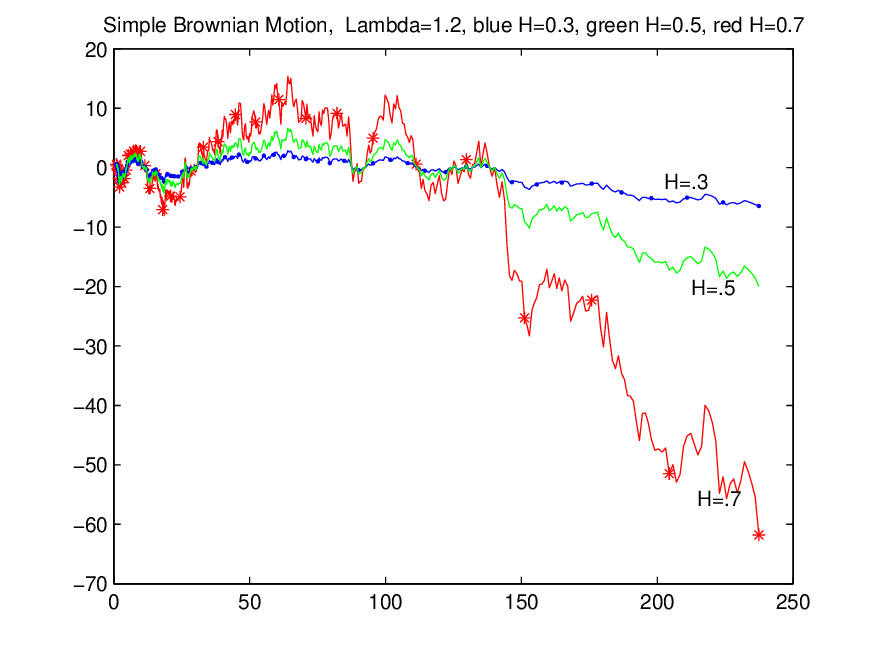}}
\caption{ \scriptsize Simple Brownian motion}
\end{figure}

\subsection{Simulation}
We present simulation of simple Brownian motion with drift which is defined in Example 3.1. For comparison it is interesting
that when the drift $g(t)$ is equal to zero, then we have simple Brownian motion with the same properties.
We have simulated and plotted simple Brownian motion and also simple Brownian motion with drift $\lambda^{n/2}\sin(\lambda^{-n}t)$
for scales $\lambda=1.2$ and $\lambda=1.5$ and different Hurst indices which have been illustrated on the Figure 1 and Figure 2.

\input{epsf}
\epsfxsize=3in \epsfysize=1.5in
\begin{figure}
\centerline{ \epsffile{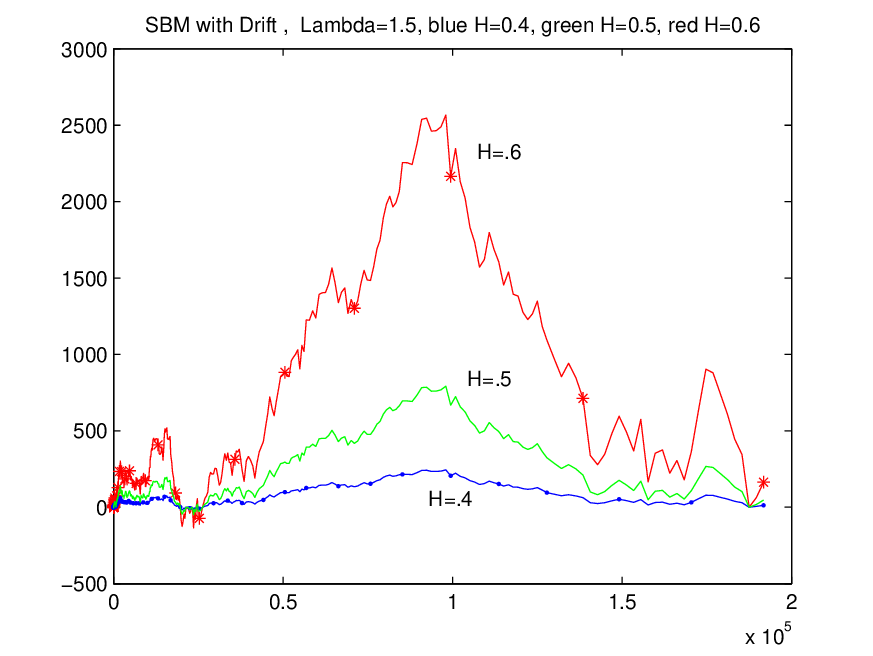}\hspace{.1in}\epsfxsize=3in \epsfysize=1.5in \epsffile{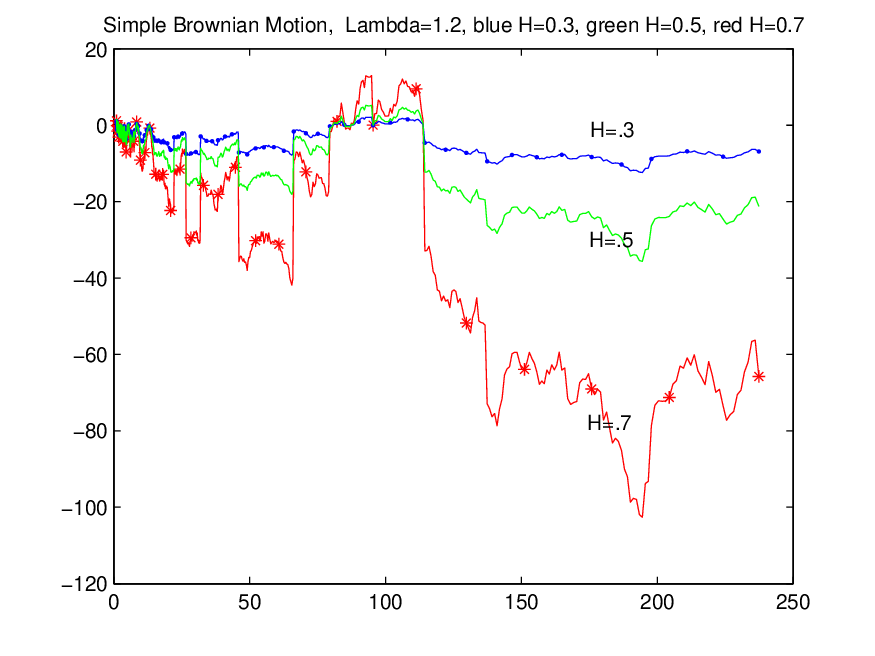}}
\caption{\scriptsize Simple Brownian motion with drift}
\end{figure}
%\vspace{1in}

\subsection{Characterization of the covariance function}
For visualizing the Theorem 3.1 in recognizing wide sense Markov property for DSI Processes, we simulated
3000 samples of the following process
$$X(t)=\sum_{n=1}^{\infty}\lambda^{n(H-\frac{1}{2})}I_{[\lambda^{n-1}, \lambda^n)}(t)B(t)$$
with Hurst index $H=0.8$ and scale $\lambda=\alpha^6$ where $\alpha=1.05$ at points $\alpha^k$, $k=0,\ldots, 3000$, where $B(t)$
is the standard Brownian motion.
So we consider samples of $M=500$ scale intervals with $T=6$ samples, by such geometric sampling, in each scale interval.
Then we estimate left hand side of (3.2) as $R_n^H(\tau)=\mbox{Cov}(X(\alpha^n), X(\alpha^{n+\tau}))$ for some different $n$ and
$\tau=kT+\nu$, say $n=9,\; k=3$ and $\nu=2 $, and the right hand side of it to examine such equality which guarantees
wide sense Markov property for such DSI process.
By applying the method, presented at the end of section 4, and by considering $3000$ samples with the above values for the parameters, we find the corresponding value for the left hand side and right hand side of (3.2) as $2.2169$ and $2.2964$ respectively, where the true value theoretically is
$2.2920$.

\subsection{Estimation of Hurst parameter}
We present a new method for estimating Hurst parameters of DSI processes.
First we need to determine the scale $\lambda$ of the process.
There are some practical methods for estimating of scale parameter in Balasis et al. \cite{bala}.
Also Rezakhah et al. \cite{ram} have considered a theoretical method in this regard, and obtained an estimation
method for scale parameter and Hurst index of semi-selfsimilar or DSI processes. They introduced a two-stage method for
estimation of Hurst index which is based on an equally spaced sampling scheme.  In their first stage the estimation was
effected by a possible self-similar behavior of the process inside each scale interval.
In this paper we consider a combination of equally spaced and geometric sampling scheme that, after determining scale intervals,
samples inside each scale interval are equally spaced but sample points in each scale interval are $\lambda$ times of the sample
points in the previous scale interval. This sampling scheme enables our estimation method that explicitly determine the Hurst
index.\\

\noindent
{\bf Estimation method}

\noindent
Our sampling scheme is to consider $T$ samples at the first scale interval $[1, \lambda)$ as equally spaced samples
$x(t_1), x(t_2), \ldots, x(t_T)$ where $ 1\leq t_1 \leq t_2 \leq \ldots \leq t_T<\lambda $ and sample points in the rest
scale intervals $(\lambda^n, \lambda^{n+1}], \; n=1, 2, \ldots, M-1; \; M\in\Bbb N$ are defined as
$x(\lambda^n t_i), \; 1\leq i \leq T$, where $\lambda$ is the scale of the DSI process. So if we denote sample points
with $\{t_i, i=1, 2, \ldots, TM\}$, then $t_{i+T}-t_i=(\lambda-1)t_i$. Also we have samples of $M$ consecutive scale intervals.
We consider first and second order variation of observations inside each scale interval as
$SS_{1,i}=\frac{1}{T-1}\sum_{k=2}^T\big(x(t_{i\cdot T+k})-x(t_{i\cdot T +k-1})\big)^2$, and $SS_{2,i}=\frac{1}{T-1}\sum_{k=3}^T\big(x(t_{i\cdot T+k})-2x(t_{i\cdot T +k-1})+x(t_{i\cdot T+k-2}\big)^2$  $i=1,2, \ldots, M-1$. Then for $j=1,2$ we evaluate
\be \mu_{j,i}=\frac{\log \left(SS_{j,i+1}/SS_{j,i}\right)}{2 \log\lambda}.\ee
Finally we provide two estimate of the Hurst index as
$\hat{H_j} = \frac{1}{M-1}\sum_{i=1}^{M-1} \mu_{j,i}$ for $j=1,2$.\\

%We follow this method by considering first order difference at lag one as
%$y(t_i)=x(t_{i+1})-x(t_i)$
%and first order difference at the scale lag as $u(t_i)=x(t_{i+T})-x(t_i)$, we also consider second order sample difference
%at lag one as $z(t_i)=x(t_{i+2})-2x(t_{i+1})-x(t_i)$ and second order difference at the scale lag as
%$v(t_i)=x(t_{i+2T})-2x(t_{i+T})-x(t_i)$. Then we find the sample variance of such different samples inside each scale interval
%and again obtain corresponding statistics as $\mu_i$ and finally evaluate their mean as the estimate of Hurst index $H$.\\

\noindent
{\bf Simulation method}

\noindent
For simulation we consider simple Brownian motion with random drift as
$$ X(t)=\sum_{i=1}^n \lambda^{i(H-\frac{1}{2})} \big(B(t)+ \lambda^{i/2} W_i\big)I_{B_i}(t),\;\;\;\; B_i=[\lambda^{i-1}, \lambda^i)$$
where $B(t)$ is the Brownian motion and $W_i$ are independent random variable with standard normal distribution. In Figure 3,
we have plotted the process and difference of order one and order two of lag one.
%at right figure, and the process with differences
%of scale lag as the left figure.
Abbreviations are, Original Observations (OO),
 %Scale Differenced Observations (SDO) and Twice Sale Differenced Observations (TSDO).
%Also
Differenced Observations (DO) and Twice Differenced Observations (TDO).
In Figure 4, 
 %we have Box-plots regarding 500 repeated estimation of different Hurst index estimation based on the above method,
%using the main samples. The left figure are Box-plots by using first order difference at scale lag of the samples. 
we have applied our estimation method and relation (5.7) for estimations of different Hurst indices by considering $50$ scale interval 
and $1000$ equally spaced samples in each scale and $100$ repetitions in each simulation and have plotted mean absolute
error (MAE) for different Hurst indices in Figure 4. 
We also applied Maximum Likelihood method for corresponding multidimensional self-similar which has the same Hurst indices as the corresponding  DSI Process and estimated Hurst indices and plotted their MAE to compare our estimation method  with this Maximum Likelihood one which based on Vidas et all [15] present good estimate for Hurst index via geometric sampling.
As it is shown in Figure 5, MAE of estimations based on  first order variation has less MAE and the MAE's of estimation by second order variations are close to it, but estimation by Maximum likelihood method has much greater MAE and are going to increase by increasing the Hurst indices.
\\ \\

\vspace{-2cm}
\input{epsf}
\epsfxsize=3in \epsfysize=1.5in
\begin{figure}\vspace{.2in}

\centerline{
%\epsffile{OBDS.ps}\hspace{.1in}\epsfxsize=3in \epsfysize=1.5in
\epsffile{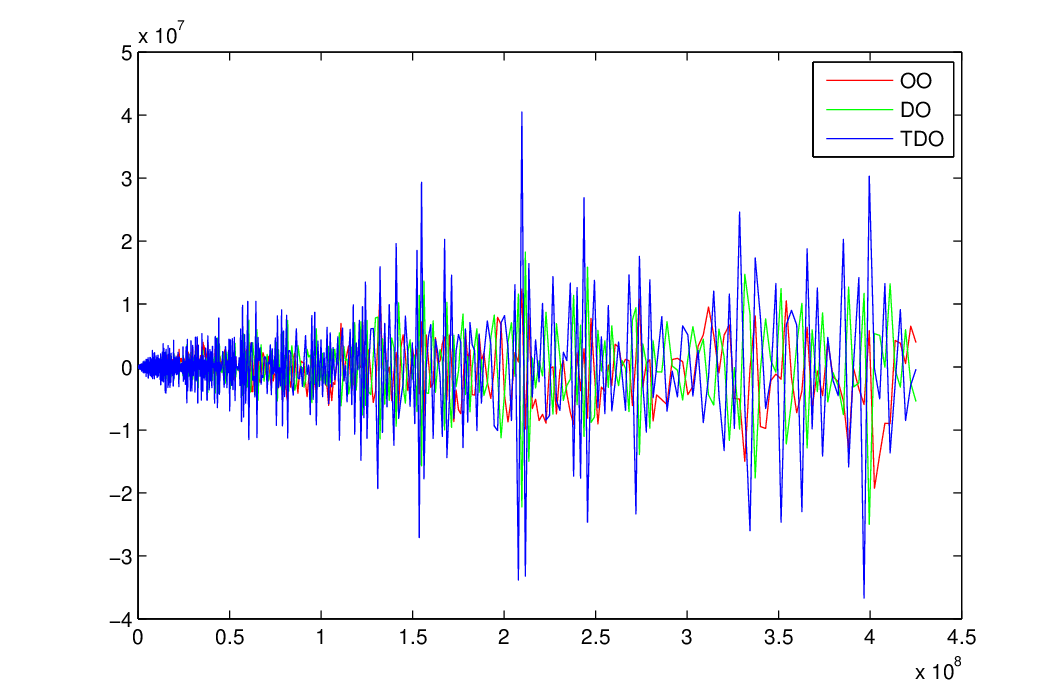}}\vspace{-.1in}

\caption{\scriptsize Simple Brownian motion with random drift. All with Hurst index H=0.8.}\vspace{.4in}

\end{figure}\vspace{.2in}

\input{epsf}
\epsfxsize=3in \epsfysize=1.5in
\begin{figure}\vspace{.2in}

\centerline{\epsffile{
%MSEDM.ps
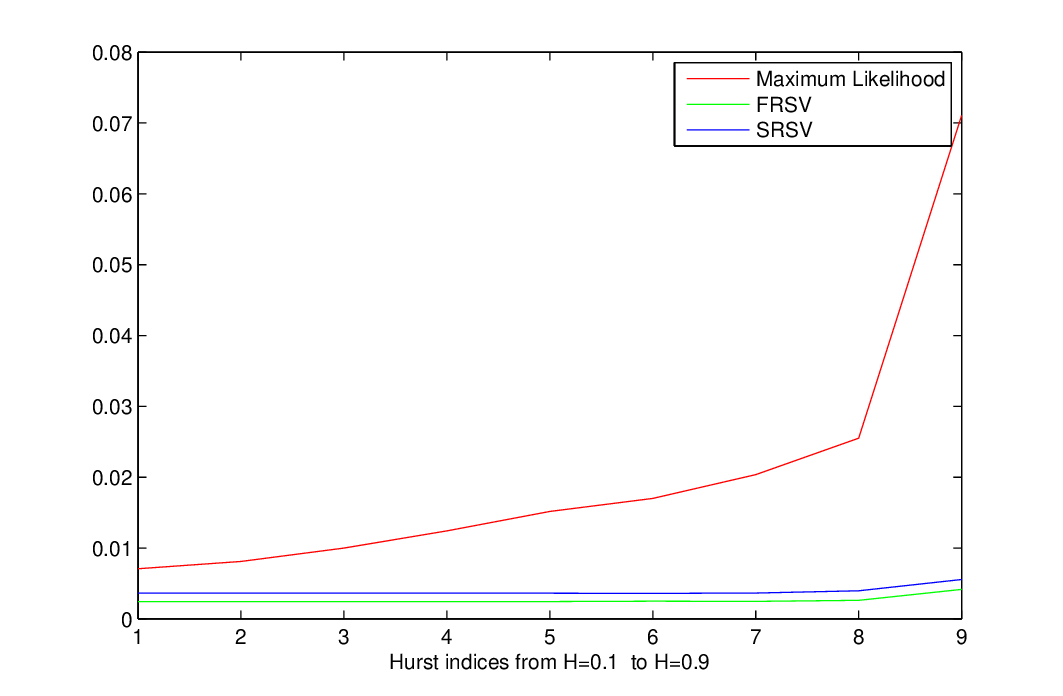}\hspace{.1in}
}\vspace{-.1in}

\caption{\scriptsize MAE for estimating Hurst index in with 100 repetitions. Horizontal axes indicates Hurst indices from 0.1 to 0.9.}
\end{figure}

\section{Conclusion}
By considering our special geometric sampling for discrete scale invariant (DSI) process we have provided a sequence of such DSI process with a correspondence multi-dimensional self-similar process.
 Imposing  Markov property for such DSI sequences, we
 have shown that the covariance function of such DSIM process is characterized by variance and covariance of adjacent  samples in the first scale interval through an explicit formula.
We have investigated this result theoretically for two examples of DSIM process as simple Brownian motion  and scale invariant autoregressive model of order one.
  By verifying this characterization formula for covariance function of simulated data we find that these estimations are very close to each other.
We have presented an efficient estimation method for the covariance function and 
spectral density matrix of corresponding multi-dimensional self-similar process as well.
Also we have proposed a new method for estimating the Hurst parameter of DSI processes by considering the ratio of sample variation of successive scale intervals. 
Comparing our estimation method for simulated data with maximum likelihood one 
  for Hurst parameter  of corresponding multi-dimensional self-similar process, reveals that our estimation method performs much better.
  This paper could initiate further research in the study of DSI process by providing such connection to the multi-dimensional processes, and also 
  has the potential to be applied for DSIM process for estimating covariance structure and spectral density matrix of corresponding multi-dimensional process.
  
%\input{epsf}
%\epsfxsize=3in \epsfysize=1.5in
%\begin{figure}\vspace{.2in}
%
%\centerline{ \epsffile{BoEDS.ps}\hspace{.1in}\epsfxsize=3in \epsfysize=1.5in\epsffile{BoESO.ps}}\vspace{-.1in}
%
%\caption{\scriptsize Box plot of different Hurst index Estimations. Left figure:  By Using Scale differenced Observations.
%Right figure: By using original observations}%\vspace{-0.2in}
%\end{figure}

\end{document}